\documentclass[12pt,a4paper,oneside]{article}
\usepackage{amssymb}
\usepackage{amsmath}
\usepackage{hyperref}
\usepackage{stmaryrd}
\usepackage{slashed}
\usepackage{url}
\usepackage{rotating}

\newcounter{count}[section]
\renewcommand{\thecount}{\arabic{section}.\arabic{count}}

\newenvironment{Umgeb1rekursiv}[1]
   {\vspace{0.5cm}
     \noindent
     \par\noindent
     \refstepcounter{count}
     \textbf{#1~\thecount}
     \hspace{0.2cm}
     \itshape
   }
   {\vspace{0.2cm}\par}

\newenvironment{Umgeb1}[1]
   {\vspace{0.5cm}
     \par\noindent
     \refstepcounter{count}
     \textbf{#1~\thecount}
     \hspace{0.2cm}
   }
   {\vspace{0.2cm}\par}

\newenvironment{UmgebBeweis}[1]
   {\vspace{0.5cm}
     \par\noindent
     \textbf{Proof.}
     \hspace{0.2cm}
   }
   {\hfill $\square$
   \vspace{0.2cm}}

\newenvironment{sat}{\begin{Umgeb1rekursiv}{Proposition}}
  {\end{Umgeb1rekursiv}}

\newenvironment{theorem}{\begin{Umgeb1rekursiv}{Theorem}}
  {\end{Umgeb1rekursiv}}

\newenvironment{bew}{\begin{UmgebBeweis}{Proof}}
   {\end{UmgebBeweis}}

\newenvironment{bem}{\begin{Umgeb1}{Remark}}
   {\end{Umgeb1}}
   
\newcommand{\R}{\mathbb{R}}

\newcommand{\C}{\mathbb{C}}
\newcommand{\N}{\mathbb{N}}
\newcommand{\Z}{\mathbb{Z}}

\newcommand{\diver}{\operatorname{div}}

\newcommand{\tr}{tr}

\newcommand{\abs}[1]{\lvert#1\rvert}
\renewcommand{\part}[2]{\frac{\partial #1}{\partial #2}}

\newcommand{\grad}{\text{grad}}

\title{On Conformal Powers of the Dirac Operator on Spin Manifolds}
\author{Matthias Fischmann}
\date{}
\begin{document}
\maketitle
\begin{abstract}
  The well known conformal covariance of the Dirac operator acting on spinor fields over a 
  semi Riemannian spin manifold 
  does not extend to powers thereof in general. For odd powers one has to 
  add lower order curvature correction terms in order to obtain conformal covariance. 
  We derive an algorithmic construction in terms of associated tractor bundles to compute 
  these correction terms. Depending on the signature of the semi Riemannian manifold in question, 
  the obtained conformal powers of the Dirac 
  operator turn out to be formally self-adjoint with respect to the $L^2-$scalar product, or formally anti-self-adjoint,
  respectively. Working out this 
  algorithm we present explicit formulas for the conformal third and fifth power of the Dirac operator. 
  
  Furthermore, we present a new family of conformally covariant differential operators acting on the 
  spin tractor bundle which 
  are induced by conformally covariant differential operators acting on the spinor bundle. 
  Finally, we will give polynomial structures for   
  the first examples of conformal powers in terms of first order differential 
  operators acting on the spinor bundle. 
\end{abstract}

\section{Introduction}

  Considering a semi Riemannian spin manifold $(M^n,g)$
  the Dirac operator is conformally covariant, see \cite{Hitchin}, whereas the Laplacian has 
  to be modified by a multiple of scalar 
  curvature, called the Yamabe operator, in order to become conformally covariant, 
  see \cite{Yamabe}, \cite{Orsted} and \cite{Branson2}. 
  Having these two examples of conformally covariant operators, Paneitz \cite{Paneitz}, actually in $1983$,  
  constructed a conformal second power of the Laplacian, i.e., he presented explicit curvature 
  correction terms for the square of 
  the Laplacian resulting in a conformally covariant operator of fourth order acting on functions. This 
  conformal second power is called the Paneitz operator. 
  Almost ten years later Graham, Jenne, 
  Mason and Sparling \cite{GJMS} constructed a series of conformally covariant 
  differential operators $P_{2N}(g)$ acting on 
  functions with leading part an $N-$th power of the Laplacian, for $N\in\N$ ($n$ odd) and $N\in\N$ 
  with $N<\frac n2$ ($n$ even). 
  The first two cases $N=1,2$ are covered by the Yamabe and the Paneitz operator. 
  Beside that construction there were two other 
  points of view describing these so-called GJMS operators. One point of view was the 
  tractor machinery used by Gover and 
  Peterson \cite{GoverPeterson} and the other one was given by Graham and Zworski 
  \cite{GZ} using a spectral theoretical point of view. 
  Again, both constructions do not produce any conformal $N-$th power of the 
  Laplacian when $n$ is even and $N>\frac n2$. 
  Although all three 
  constructions are algorithmic explicit formulas have very rarely been produced, 
  due to their complexity. 
  In case of Einstein manifolds, Gover \cite{Gover} proved a product structure of 
  shifted Laplacains of the GJMS operators. 
  Recent results of Juhl \cite{Juhl2, Juhl1} simplified the structure by showing that 
  the GJMS operators can be described as polynomials in second order differential operators. 

  Let us now move to the spinor case: It follows from \cite[Theorem $8.13$]{Slovak} that no conformal even 
  powers of the Dirac operator can be expected. Holland and Sparling 
  \cite{HS} proved the existence of conformal odd powers of the Dirac operator. 
  In the even dimensional case, their 
  construction failed to give conformal odd powers when the order exceed the dimension. 
  The first explicit formula for a conformal third power is due to Branson \cite{BransonB}, 
  which he derived using 
  tractor techniques. Later on, Gillarmou, Moroianu and Park \cite{GMP1} 
  gave a construction for 
  conformal odd powers of the Dirac operator using a spectral theoretical point of view. However, in the 
  even dimensional case, this 
  does not yield conformal powers when the order exceed the dimension. 
  They also gave an explicit formula for the 
  conformal third power of the Dirac operator, in agreement with the result of Branson. In \cite{ES}, 
  Eelbode and 
  Sou\v{c}ek derived a product structure of shifted Dirac operators for conformal powers of the Dirac operator in 
  case of the Riemannian sphere. But in general, due to the complexity of the 
  underlying algorithms, further examples were not known in the literature. 

  The mentioned constructions of conformal powers of the Laplacian and the Dirac operator based on the 
  ambient metric construction, introduced by Fefferman and Graham \cite{FG, FG3}. In general, the construction 
  of the ambient metric is obstructed in case of even dimensional manifolds. This is the reason that in 
  those dimensions the conformal powers of the Laplacian and Dirac operator only exist up to 
  the order mentioned above.

  \vspace{5pt}
  The paper is organized as follows. We always assume that $(M,g)$ is a semi Riemannian spin manifold. 
  
  In Section $2$ we recall basic notation from semi Riemannian geometry and spin geometry. Furthermore, 
  we recall parabolic geometries with main focus conformal geometry. That means, we will present the 
  standard tractor bundle with is normal conformal Cartan connection. 
  This construction goes back to Cartan \cite{Cartan} and 
  Thomas \cite{Thomas} and was put into a modern language by \v{C}ap and Slov\'{a}k \cite{CapSlovak}. 
  Dealing with conformal spin structures naturally leads to the spin tractor bundle, which is also introduced. 
  
  In Section $3$ we recall 
  the construction of so-called splitting operators, using Casimir techniques \cite{CapSoucek}. They will 
  be used for 
  the construction of a series of conformally covariant differential operators $P_{2N}^{\mathcal{S}(M)}(g)$ 
  acting on the spin tractor bundle 
  by translation of the strongly invariant Yamabe operator in the sence of \cite{EastwoodRice}.

  In Section $4$ we use the splitting operators to construct 
  conformal odd powers of the Dirac operator, again using the curved translation principle of Eastwood and Rice. 
  In case of even dimensional manifolds this construction does not give any conformal odd 
  powers when the order exceeds the dimension. Furthermore, depending on the signature of metric, 
  we prove that the constructed operators are formally self-adjoint, or anti-self-adjoint, with 
  respect to the $L^2-$scalar product, respectively. In the special case 
  of Einstein manifolds, we prove that the first examples of conformal powers of the Dirac operator posseses a 
  product structure, consisting of shifted Dirac operators. We then return to the general setting, and show 
  that the splitting operators 
  can be used to construct a new family of conformally covariant differential 
  operators $L_k(g)$, for $k\in 2\N+1$, acting on the spin tractor bundle. These differ sligthly 
  from the $P_{2N}^{\mathcal{S}(M)}(g)$, however they do have the same conformal bi-degree. 
  Finally, we give a new polynomial structure for the first examples of the conformal 
  powers of the Dirac operator, analogous to the work of Juhl in case of the GJMS operators. 
  For computations which are omitted and further references we refer to \cite{Fischmann}.

  \textbf{Acknowledgements:} I like to thank the BMS, SFB $647$ and Eduard \v{C}ech Institute 
  for their financial support. 
  Furthermore, I would like to take this opportunity to express my gratitude to Helga Baum and 
  Andreas Juhl. Helga Baum 
  introduced me to the realm of semi Riemannian geometry, especially conformal geometry, 
  whereas Andreas Juhl inspired 
  me to work on the subject of conformally 
  covariant differential operators.

\section{Preliminaries}

  Let $(M,g)$ be a semi Riemannian spin manifold of signature $(p,q)$. We begin by fixing 
  some curvature conventions and introduce tensors fields 
  which will be used throughout the paper. Next, we recall the concept of spinor bundles associated to 
  $(M,g)$. A detailed treatment of spinor bundles and tools used within 
  the paper can be found in \cite{LM, Baum1}. We then go on to recall the concept of conformally 
  covariant differential operators
  in the sence of \cite{Kosmann}. Finally, we present a conformal invariant calculus in the language of parabolic 
  geometry, see \cite{CapSlovakSoucek1,CapSlovakSoucek} and \cite{CapSlovak}, 
  and related tractor bundles, upon which our construction of conformal powers of the Dirac operator is based. 

  \subsection{Tensor conventions}
    Let us denote by $\nabla^{LC}:\Gamma(TM)\to\Gamma(T^*M\otimes TM)$ the Levi-Civita connection 
    canonically associated to $(M,g)$. The curvature tensor of 
    the Levi-Civita connection is defined by 
    $R(X,Y)Z:=\nabla^{LC}_X\nabla^{LC}_YZ-\nabla^{LC}_Y\nabla^{LC}_XZ-\nabla^{LC}_{[X,Y]}Z$, and 
    the Riemannian curvature tensor is defined by 
    $\mathcal{R}(X,Y,Z,W):=g(R(X,Y)Z,W)$, for $X,Y,Z,W\in\mathfrak{X}(M)$. 
    Further tensor fields  which can be built from  
    the Riemannian curvature tensor (using covariant derivatives and contractions) are:
    \begin{itemize}
      \item $Ric(X,Y):=\tr_g\left(\mathcal{R}(X,\cdot,\cdot,Y)\right)$ (Ricci tensor),
      \item $\tau:=\tr_g\left(Ric(\cdot,\cdot)\right)$ (scalar curvature),
      \item $J:=\frac{1}{2(n-1)}\tau$ (normalized scalar curvature),
      \item $P(X,Y):=\frac{1}{n-2}\left(Ric(X,Y)-Jg(X,Y)\right)$ (Schouten tensor),
      \item $W(X,Y,Z,W):=\mathcal{R}(X,Y,Z,W)+P\owedge g(X,Y,Z,W)$ (Weyl tensor),
      \item $C(X,Y,Z):=\nabla^{LC}_XP(Y,Z)-\nabla^{LC}_YP(X,Z)$ (Cotton tensor),
      \item $B(X,Y):=\tr_g\left(\nabla^{LC}_\cdot C(\cdot,X,Y)\right)+g\left(P(\cdot,\cdot),W(\cdot,X,Y,\cdot)\right)$ 
        (Bach tensor),
    \end{itemize}
    where the Kulkarni-Nomizu product $\owedge$ is defined by 
    \begin{align*}
      P\owedge g(X,Y,Z,W):=&P(X,Z)g(Y,W)+P(Y,W)g(X,Z)\\
      &-P(X,W)g(Y,Z)-P(Y,Z)g(X,W),
    \end{align*}
    for $X,Y,Z,W\in\mathfrak{X}(M)$. Finally, the semi Riemannian metric yields the usual isomorphisms 
    $\cdot^\natural:T^*M\to TM$ and $\cdot^\flat: TM\to T^*M$.

  \subsection{Clifford algebras, spin groups and their  representations}
    Consider the vector space $\R^n$ ($n=p+q$) together with the scalar product 
    $\langle\cdot,\cdot\rangle_{p,q}$ of index $p$, i.e., $\langle e_i,e_j\rangle_{p,q}=\varepsilon_i\delta_{ij}$, 
    where $\{e_i\}$ 
    is the standard basis 
    of $\R^n$, $\varepsilon_i=-1$ for $1\leq i\leq p$; $\varepsilon_i=1$, for $p+1\leq i\leq n$, 
    and $\delta_{ij}$ denotes the Kronecker delta. 
    Consider the Clifford algebra of $\R^{p,q}:=(\R^n,\langle\cdot,\cdot\rangle_{p,q})$ realized 
    by $\mathcal{C}_{p,q}:=T(\R^n) /J$,
    where $T(\R^n)$ denotes the tensor algebra of $\R^n$, and $J$ is the two-sided ideal in 
    $T(\R^n)$ generated by the relations $x\otimes x=-\langle x,x\rangle_{p,q}$, for $x\in\R^n$. 
    The Clifford algebra carries a $\Z_2-$grading, given by even 
    and odd elements, i.e., $\mathcal{C}_{p,q}=\mathcal{C}^0_{p,q}\oplus\mathcal{C}^1_{p,q}$. 
    We denote the group of units of $\mathcal{C}_{p,q}$ by $\mathcal{C}_{p,q}^*$ and call it the Clifford group. 
    This leads to two important subgroups, 
    the pin group $Pin(p,q)$, given by products of elements $x\in\R^n$
    of lenght $\pm 1$, and the spin group $Spin(p,q):=Pin(p,q)\cap\mathcal{C}^0_{p,q}$. There is an 
    algebra isomorphism of the comlexified Clifford algebra 
    \begin{align*}
     \Phi^{even/odd}_{p,q}:\mathcal{C}_{p,q}^{\C}\to\begin{cases}
                                    Mat(2^m,\C),\quad &n=2m\\
                                    Mat(2^m,\C)\oplus Mat(2^m,\C),\quad&n=2m+1
                                  \end{cases}.
    \end{align*}
    It is defined as follows: Set
    \begin{align*}
      g_1:=\begin{pmatrix} i&0\\0&-i\end{pmatrix},\quad g_2:=\begin{pmatrix} 0&i\\i&0\end{pmatrix},\quad 
      T:=\begin{pmatrix}0&-i\\i&0\end{pmatrix},\quad E:=\begin{pmatrix} 1&0\\0&1\end{pmatrix}
    \end{align*}
    and 
    \begin{align*}
      \alpha(j):=\begin{cases} 1,\quad j\in 2\N-1\\ 2,\quad j\in 2\N\end{cases},\quad
      \tau(j):=\begin{cases} i,\quad j\leq p\\ 1,\quad j>p\end{cases}.
    \end{align*}
    In the case of $n=2m$, we use an orthonormal basis $\{e_i\}$ of $\R^{p,q}$ to define 
    the isomorphism 
    \begin{align*}
      \Phi^{even}_{p,q}(e_j):=\tau(j) E\otimes\ldots\otimes E\otimes g_{\alpha(j)}\otimes T\ldots\otimes T.
    \end{align*}
    Here, the right hand side is a product of $m$ matrices, $\left[\frac{j-1}{2}\right]$ of them are copies of $T$, 
    and the tensor product used is the Kronecker tensor product for matrices. 
    In the case of $n=2m+1$, we set
    \begin{align*}
      \Phi^{odd}_{p,q}(e_j)
         :=&\left(\Phi^{even}_{p,q-1}(e_j),\Phi^{even}_{p,q-1}(e_j)\right),\quad j=1,\ldots 2m,\\
      \Phi^{odd}_{p,q}(e_{2m+1}):=&(iT\otimes\ldots\otimes T,-iT\otimes\ldots\otimes T),
    \end{align*}
    where $\{e_i\}$ is an orthonormal basis of $\R^{p,q}$. 
    Hence, in the case of $n=2m$, the Clifford algebra $\mathcal{C}^{\C}_{p,q}$ 
    has (up to equivalence) an unique irreducible 
    representation $\Phi_{p,q}:=\Phi^{even}_{p,q}$, whereas, in the case of $n=2m+1$ 
    it has (up to equivalence) two unique irreducible representations denoted by $\Phi^0_{p,q}$ 
    and $\Phi^1_{p,q}$. 
    In all cases the representation space is $\Delta_{p,q}:=\C^{2^m}$. Note that in the 
    case of $n=2m+1$, both irreducible 
    representations $\Phi^0_{p,q}$ and $\Phi_{p,q}^1$ become equivalent when they are restricted 
    to the even part $\mathcal{C}_{p,q}^0$. 
    Restricting $\Phi_{p,q}$, in the even case, or $\Phi_{p,q}^0$, in the odd case, 
    to the spin group yields a representation of the spin group, which 
    will be denoted by $\kappa_{p,q}$. 
    This is the spinor representation we will work with. 
    Again, in the case of $n=2m$ we have that  $\kappa_{p,q}$ decomposes into two non-equivalent irreducible 
    representations, whereas in the case of $n=2m+1$ the representation 
    $\kappa_{p,q}$ is irreducible. 

    On the representation space $\Delta_{p,q}$ there exists a $Spin_0(p,q)-$invariant hermitian scalar product 
    $(v,w)_\Delta:=(b\cdot v,w)_{\C^{2^m}}$, where $Spin_0(p,q)$ denotes the connected component 
    containing the identity, 
    $(\cdot ,\cdot)_{\C^{2^m}}$ is the standard 
    hermitian scalar product on $\C^{2^m}$, and 
    $b:=i^{\frac{p(p-1)}{2}}e_1\cdot\ldots\cdot e_p$. For Riemannian 
    signature (that is $p=0$) it reduces to the standard hermitian scalar product, which 
    is $Spin_0(0,n)=Spin(0,n)-$invariant.

  \subsection{Spin structures and spinor bundles}\label{Preliminaries}
    Let $(\mathcal{Q}^g,f^g)$ be a spin structure for $(M,g)$, i.e., a $\lambda-$reduction of the orthonormal frame 
    bundle $(\mathcal{P}^g,\pi,M,SO(p,q))$, where $\lambda:Spin(p,q)\to SO(p,q)$ denotes the usual 
    twofold covering of $SO(p,q)$. The associated vector bundle 
    $S(M,g):=\mathcal{Q}^g\times_{(Spin_0(p,q),\kappa_{p,q})}\Delta_{p,q}$ over $M$ is called the spinor 
    bundle of $(M,g)$. The hermitian scalar product $(\cdot,\cdot)_\Delta$ induces a scalar product on the spinor 
    bundle by $<\psi,\phi>:=(v,u)_{\Delta}$, for $\psi=[q,v], \phi=[q,u]\in S(M,g)$. 
    Due to the reduction property of $(\mathcal{Q}^g,f^g)$ we obtain an isomorphism 
    $TM\simeq \mathcal{Q}^g\times_{(Spin_0(p,q),\rho\circ\lambda)}\R^n$, 
    where $\rho$ denotes the standard representation 
    of $SO(p,q)$ on $\R^{p,q}$, and thus we may define the Clifford 
    multiplication $\mu:TM\otimes S(M,g)\to S(M,g)$ by 
    \begin{align*}
      \mu(X\otimes\psi):=\begin{cases}[q,\Phi_{p,q}(x)v],\quad &n=2m\\ 
        [q,\Phi^0_{p,q}(x)v],\quad & n=2m+1, \end{cases}
    \end{align*}
    for $X=[q,x]\in TM$, and $\psi=[q,v]\in S(M,g)$. If there is no confusion we will use $X\cdot\psi$ 
    instead of $\mu(X\otimes\psi)$. Clifford multiplication 
    extends to the exterior algebra of $T^*M$ by 
    \begin{align*}
      w\cdot\psi:=\sum_{1\leq i_1<\ldots<i_k\leq n}\varepsilon_{i_1}\ldots\varepsilon_{i_k}
         w(s_{i_1},\ldots,s_{i_k})s_{i_1}\cdot\ldots\cdot s_{i_k}\cdot\psi,
    \end{align*}
    where $w\in(\Lambda^kM)_x$, $\psi\in S(M,g)_x$ and $\{s_i\}$ is an orthonormal basis in $T_xM$, for $x$ 
    the base point. Note that 
    Clifford multiplication varies smoothly on $M$, thus it descents to sections of corresponding vector 
    bundles. In order to define a covariant derivative 
    on the spinor bundle in a canonical way we choose the Levi-Civita connection form 
    $A^g\in\Omega^1\left(\mathcal{P}^g,\mathfrak{so}(p,q)\right)$, induced by $\nabla^{LC}$, 
    and define, using the isomorphism $\lambda_*:\mathfrak{spin}(p,q)\to\mathfrak{so}(p,q)$ 
    (the differential of the covering map at the identity), a connection form
    $\tilde{A}^g:=\lambda_*^{-1}\circ A^g\circ \lambda_*\in\Omega^1\left(\mathcal{Q}^g,\mathfrak{spin}(p,q)\right)$ 
    on $\mathcal{Q}^g$. This induces a covariant derivative on the associated 
    vector bundle $S(M,g)$ in the usual way, i.e., locally we have
    \begin{align*}
      \nabla^{S(M,g)}_X\psi_{\mid U}=&\left[q,dv_q(X)+(\kappa_{p,q})_*\left((\tilde{A}^g)^q(X)\right)v\right]\\
      \overset{loc.}{=}&X(\psi)+\frac 12\sum_{i<j}\varepsilon_i\varepsilon_jg(\nabla^{LC}_Xs_i,s_j)s_i\cdot s_j\cdot\psi,
    \end{align*}
    for local sections $\psi=[q,v]:U\to S(M,g)$ and $s=\{s_i\}:U\to \mathcal{P}^g$, and 
    for $X(\psi):=[q,dv_q(X)]$. Here, 
    $(\tilde{A}^g)^q(X):=\tilde{A}^g_q(dq(X))$ is the connection $1-$form induced by the 
    local section  $q:U\to \mathcal{Q}^g$. 
    The covariant derivative 
    $\nabla^{S(M,g)}$ leads to the definition of the Dirac operator by 
    \begin{align*}
     \slashed{D}:\Gamma\left(S(M,g)\right)&\to\Gamma\left(S(M,g)\right)\\
     \psi&\mapsto \slashed{D}\psi:=\mu\left((\nabla^{S(M,g)}\psi)^\natural\right),
    \end{align*}
    where $\cdot^\natural$ indicates the identification $T^*M\simeq TM$ induced by $g$. 
    Locally the Dirac operator reads 
    $\slashed{D}\psi\overset{loc.}{=}\sum_i\varepsilon_i s_i\cdot\nabla^{S(M,g)}_{s_i}\psi$. The following 
    list collects 
    useful formulas, some are well known, see \cite{Baum1,LM}, and while the remainder 
    are straightforward to derive:
    For $\psi,\phi\in\Gamma\left(S(M,g)\right)$ and $X,Y\in\mathfrak{X}(M)$, one has
    \begin{itemize}
        \item[(1)] $\nabla^{S(M,g)}_X(Y\cdot\psi)=\nabla^{LC}_XY\cdot\psi+Y\cdot\nabla^{S(M,g)}_X\psi$,
        \item[(2)] $\nabla^{S(M,g)}$ is metric with respect to $<\cdot,\cdot>$,
        \item[(3)] $<X\cdot\psi,\phi>+(-1)^p<\psi,X\cdot\phi>=0$,
        \item[(4)] $\mathcal{R}^{S(M,g)}(X,Y)\psi=\frac 12 \mathcal{R}(X,Y)\cdot\psi$, where the 
          Riemannian curvature tensor is considered as endomorphism of $2-$forms,
        \item[(5)] $[\slashed{D},f]\psi=\slashed{D}(f\psi)-f\slashed{D}\psi=\grad^g f\cdot\psi$, 
          for any $f\in\mathcal{C}{^\infty}(M)$,
        \item[(6)] $[\slashed{D},\nabla^{S(M,g)}]\psi=\frac 12 Ric(X)^\natural\cdot\psi$ and 
        \item[(7)] $\slashed{D}^2\psi=-\Delta^{S(M,g)}_g\psi+\frac{\tau}{4}\psi$ is the Bochner formula, where 
          $\Delta^{S(M,g)}_g:=\tr_g(\nabla^{T^*M\otimes S(M,g)}\circ \nabla^{S(M,g)})$ is the Bochner 
          Laplacian on spinor fields. 
    \end{itemize}

    Concerning questions of self-adjointness of certain operators on spinor fields 
    we introduce a bracket notation. Let $T$ be a symmetric $(0,2)-$tensor and $\psi$ a spinor field. 
    We define first a $1-$form $T\cdot\psi$ with values in the spinor 
    bundle by $T\cdot\psi(X):=T(X)^\natural\cdot\psi$. 
    Then the following brackets are defined:
    \begin{align}\label{eq:bracket1}
      (T,\nabla\psi):=&\mu\left(\tr_g(T(\cdot)^\natural\otimes\nabla_\cdot\psi))\right)
         \overset{loc.}{=}\sum_i\varepsilon_i T(s_i)^\natural\cdot\nabla_{s_i}\psi,\\
      (\nabla,T\cdot\psi):=&-\delta^{\nabla^{S(M,g)}}(T\cdot\psi),
    \end{align}
    where, for $\eta\in\Omega^1(M,S(M,g))$, 
    $\delta^{\nabla^{S(M,g)}}\eta\overset{loc.}{=}-\sum_i\varepsilon_i(\nabla^{S(M,g)}_{s_i}\eta)(s_i)$ is the 
    co-differential of $d^{\nabla^{S(M,g)}}$. Note that the last 
    bracket can be rewritten as 
    \begin{align*}
      (\nabla^{S(M,g)},T\cdot\psi)=(T,\nabla^{S(M,g)}\psi)-(\delta^{\nabla^{LC}}T^\natural)\cdot\psi,
    \end{align*}
    where $\delta^{\nabla^{LC}}$ denotes the co-differential of $d^{\nabla^{LC}}$. 
    Next, we define a $(0,2)-$tensor $T^2$ by 
    $T^2(X,Y):=T(T(X)^\natural,Y)$, and a further bracket by 
    \begin{align}
      (C,P\cdot\psi):=&\sum_i\varepsilon_i C(s_i)\cdot P(s_i)\cdot\psi,\label{eq:bracket3}
    \end{align}
    where the Cotton tensor is considered as $C(X):=C(\cdot,\cdot,X)\in\Omega^2(M)$. 
    Analogously one defines $(P,C\cdot\psi)$. Using the same notation for those brackets will not 
    lead to any confusion. 
    Two more product types, needed later on, are 
    \begin{align}
      W\cdot W\cdot\psi:=&\sum_{i,j}\varepsilon_i\varepsilon_j W(s_i,s_j)\cdot W(s_i,s_j)\cdot\psi,\\
      C\cdot W\cdot\psi:=&\sum_{i,j}\varepsilon_i\varepsilon_j C(s_i,s_j,\cdot)^\natural\cdot W(s_i,s_j)\cdot\psi,
    \end{align}
    where Clifford multiplication of $2-$forms $W(X,Y):=W(X,Y,\cdot,\cdot)\in\Omega^2(M)$ appears. 
    Similary we define $W\cdot C\cdot\psi$.

  \subsection{Conformal structures and conformally covariant differential operators}
    We say that another metric $\hat{g}$ on $M$ is conformally 
    related to $g$ if there is a smooth function 
    $\sigma\in\mathcal{C}^\infty(M)$ such that $\hat{g}=e^{2\sigma}g$. This clearly defines an equivalence relation 
    among metrics on $M$. 
    We call $(M,c:=[g])$ a conformal semi Riemannian manifold. 
    Note that signature and orientation are invariant under 
    a conformal change of a metric. A conformal structure $[g]$ on $M$ induces a 
    $CO(p,q)\simeq\R^+\times SO(p,q)-$reduction $(\mathcal{P}^0,\pi,M,CO(p,q))$ of 
    the frame bundle $(GL(M),\pi,M,Gl(n,\R))$, in analogly to semi Riemannian structures $g$ on $M$ where 
    $GL(M)$ reduces to the orthonormal frame bundle $\mathcal{P}^g$. We should point out, that in contrast to 
    the semi Riemannian case there is no distinguished connection form on the conformal frame bundle, but there 
    is one on its first prolongation which will be discussed in the next subsection. 
    
   We will now define a conformal spin structure on a conformal manifold $(M,c)$. 
   Consider the conformal spin group
    \begin{align*}
      CSpin(p,q):=\R^+\times Spin(p,q)
    \end{align*}
    and the map $\lambda^c:CSpin(p,q)\to CO(p,q)$, defined by $\lambda^c(a,g):=a\lambda(g)$. 
    A conformal spin structure $(\mathcal{Q}^0,f^0)$ on $(M,c)$ is defined to be 
    a $\lambda^c-$reduction of 
    the conformal frame bundle. Conformal spin structures on $(M,c)$ are equivalent to spin structures on 
    $(M,g)$ in the following way: Given a spin structure $(\mathcal{Q}^g,f^g)$ on $(M,g)$ we 
    define a conformal spin structure $(\mathcal{Q}^0, f^0)$ on $(M,c)$ by taking the extension 
    $\mathcal{Q}^0:=\mathcal{Q}^g\times_{Spin(p,q)}CSpin(p,q)$, and setting $f^0:=f^g\times\lambda^c$. 
    Conversely, given a conformal spin structure $(\mathcal{Q}^0,f^0)$ 
    on $(M,c)$, choosing $g\in c$, 
    we define, using the obvious reduction map $\iota:\mathcal{P}^g\to \mathcal{P}^0$, a 
    spin structure $(\mathcal{Q}^g,f^g)$ on $(M,g)$ by 
    $\mathcal{Q}^g:=\{q\in\mathcal{Q}^0\mid f^0(q)\in\iota(\mathcal{P}^g)\}$ and 
    $f^g:=f^0_{\mid\mathcal{Q}^g}$.
    
    \begin{bem}
     Since we have no distinguished connection form on the conformal frame bundle we cannot 
     build up a conformally invariant differential calulus 
     on the tangent bundle. However, as we will see in the next two subsections, 
     there is a first prolongation of the conformal frame 
     bundles which possess a distinguished Cartan connection. 
    This Cartan connection induces a covariant derivative on the so-called 
     tractor bundles. Then, by fixing a representative $g\in c$, it is possible to identify within that 
     covariant derivative, its curvature, or in the divergence of its curvature tensors like Schouten, Weyl, 
     Cotton and Bach associated to $g$.
    \end{bem}

    Let us finish this subsection with the notion of conformally covariant differential operators 
    acting between sections of two  
    vector bundles $E\to M$ and $F\to M$ over $(M,g)$. We say that a linear differential 
    operator $D(g):\Gamma(E)\to\Gamma(F)$ 
    is $g-$geometrical if it is a polynomial in $g$, $g^{-1}$, $\nabla^{LC}$ and $\mathcal{R}$. A 
    $g-$geometrical differential 
    operator $D(g)$ is said to be conformally covariant of bi-degree $(a,b)$ if there exists $a,b\in \R$ such that
    \begin{align*}
      D(e^{2\sigma}g)(e^{a\sigma}\psi)=e^{b\sigma}D(g)\psi,
    \end{align*}
    for any metric $e^{2\sigma}g$, and $\psi\in\Gamma(E)$. If the bundles $E$ 
    and $F$ depend on the chosen metric, but can be related by a bundle map for conformally related metrics, 
    then this map can be used to define conformally covariant operators between $E$ and $F$. 
    An example is given by the spinor bundle $S(M,g)$; here,
    there exists a bundle isomorphism $F_\sigma:S(M,g)\to S(M,e^{2\sigma}g)$ induced from the map 
    $\Lambda_\sigma:\mathcal{P}^g\to\mathcal{P}^{e^{2\sigma}g}$ (which is given by 
    $\Lambda_\sigma(s_1,\ldots,s_n):=(e^{-\sigma}s_1,\ldots,e^{-\sigma}s_n)$), 
    and the covering property of spin structures, see \cite{Baum1}. Another example is given by 
    the maps $T(g,\sigma)$ and $T^{\mathcal{S}(M)}(g,\sigma)$, see Subsection \ref{tractorbundles}, these 
    identify the metric decomposition of certain tractor bundles with respect to two 
    representatives from the conformal class. 
    Examples of conformally covariant operators are the Yamabe operator acting on functions, the Dirac operator 
    and the twistor operator acting on spinor fields. 
    In Section \ref{relevantoperators} and \ref{results} we will deal with more conformally 
    covariant differential operators.

  \subsection{Parabolic geometries for conformal spin structures}

    Parabolic geometries are special classes of Cartan geometries, which themselfs are 
    curved versions of Klein geometries $(G,\pi,G/H,H;w_G)$, where 
    $G$ is a Lie group, $H\subset G$ is a closed subgroup, and $w_G$ is the Maurer-Cartan form. 
    
    For $H\subset G$ as above and $M$ a smooth manifold, a 
    Cartan geometry $(\mathcal{G},\pi,M,H;w)$ of type $(G,H)$, consists of 
    an $H-$principal bundle $\mathcal{G}$ over $M$ with a Cartan connection 
    $w\in\Omega^1(\mathcal{G},\mathfrak{g})$, such that $(1)$ $w(\tilde{X})=X$ for every 
    $X\in \mathfrak{h}$ (where $\tilde{X}$ denotes the fundamental vector field of $X$), $(2)$ 
    $w:T_u\mathcal{G}\to\mathfrak{g}$ is an isomorphism, for every $u\in\mathcal{G}$, and 
    $(3)$ $(R_h)^*w=Ad(h^{-1})\circ w$, for every $h\in H$. 

    A Cartan geometry $(\mathcal{G},w)$ of type $(G,H)$, for which $H$ is a parabolic 
    subgroup inside a semisimple Lie group $G$, is referred to 
    as a parabolic geometry. For more details see \cite{Sharpe} and \cite{CapSlovak}. 

    A conformal manifold $(M,c)$ of signature $(p,q)$ can be described as a parabolic geometry 
    as follows: Let us denote 
    $G:=O(p+1,q+1)/\{\pm Id\}$ the projective orthonormal group. 
    In terms of the standard orthonormal basis $\{e_\alpha\}_{\alpha=0}^{n+1}$ with respect to 
    the standard semi Riemannian metric $\langle\cdot,\cdot\rangle_{p+1,q+1}$ on $\R^{n+2}$, we define 
    the following basis 
    \begin{align*}
      f_0:=\frac{1}{\sqrt{2}}(e_{n+1}-e_0),\quad f_i:=e_i,\quad f_{n+1}:=\frac{1}{\sqrt{2}}(e_{n+1}+e_0)
    \end{align*} 
    on $\R^{n+2}$. The stabilizer $B:=stab_{\R f_0}(G)$ of the isotropic line $\R f_0$ defines a parabolic 
    subgroup of $G$, and it is isomorphic, under the projection $O(p+1,q+1)\to G$, 
    to the following subgroup of $O(p+1,q+1)$:
    \begin{align*}
       B\simeq\left.\left\{ Z(a,A,v):=\begin{pmatrix}
          a^{-1}&v^{t}&b\\
          0&A&x\\
          0&0&a
        \end{pmatrix}\right|
        \begin{matrix}
          a\in\R^{+},v\in\R^{p,q},A\in O(p,q),\\
          x:=-aAJ^{p,q}v,\\
          b:=-\frac 12 a\langle v,v\rangle_{p,q}
        \end{matrix}
       \right\},
    \end{align*} 
    where $J^{p,q}:=diag(-I_p,I_q)$ and $I_r$ denotes the identity matrix of size $r$. This group carries a 
    semi direct product structure: $B\simeq B_0\ltimes_\rho B_1$ for 
    \begin{align*}
       B_{0}:=&\{X(a,A):=Z(a,A,0)\in B\}\simeq CO(p,q),\\
      B_{1}:=&\{Y(v):=Z(1,I_{n},v)\in B\}\simeq \R^{n},
    \end{align*}
    where $\rho: B_0\times B_1\to B_1$ is the conjugation map $\rho(b_0)b_1:=b_0b_1b_0^{-1}$. Finally, let us 
    denote $B_{-1}:=\{Y(v)^t\mid v\in \R^n\}$. This will be needed for the grading of the Lie algebra of 
    $G$, i.e., $\mathfrak{g}:=LA(G)=\mathfrak{b}_{-1}\oplus \mathfrak{b}_0 \oplus\mathfrak{b}_1$ is 
    a $\abs{1}-$graded Lie algebra. 
    In terms of matrices one has 
    \begin{align*}
      \mathfrak{g}=\left.\left\{M(x,(A,a),z):=\begin{pmatrix}
                                                -a&z&0\\
                                                x&A&-J^{p,q}z^t\\
                                                0&-J^{p,q}x^t&a  
                                              \end{pmatrix}\right| 
        \begin{matrix} x\in\R^n,z\in(\R^n)^*, \\ a\in\R, A\in\mathfrak{o}(p,q)\end{matrix} \right\}
    \end{align*}
    and 
    \begin{align*}
      \mathfrak{b}_{-1}&=\{M(x,(0,0),0)\in\mathfrak{g}\}\simeq \R^n,\\
      \mathfrak{b}_0&=\{M(0,(a,A),0)\in\mathfrak{g}\}\simeq \mathfrak{co}(p,q),\\
      \mathfrak{b}_{1}&=\{M(0,(0,0),z)\in\mathfrak{g}\}\simeq (\R^n)^*.
    \end{align*}
    In this setting it is shown in \cite[Section $1.6$]{CapSlovak} that there exists a parabolic 
    geometry $(\mathcal{P}^1,w^{nc})$ of type $(G,B)$ uniquely associated to the conformal structure. 
    Roughly speeking, 
    the $B-$principal bundle $\mathcal{P}^1$, called the first prolongation of the conformal frame bundle, 
    is the collection of horizontal and torsion free subspaces in $T\mathcal{P}^0$, 
    and the normal conformal Cartan connection 
    $w^{nc}$ is an extension of the soldering form of $\mathcal{P}^1$.    
    Additionally, one has that $(\mathcal{P}^1,\pi^1,\mathcal{P}^0,B_1)$ is a $B_1-$principal 
    bundle over $\mathcal{P}^0$, whereas 
    $(\mathcal{P}^1,\pi^0,M,B)$ is a $B-$principal bundle over $M$, with the obvious projection maps. 

    As we promised earlier, choosing a metric $g$ from the conformal class, we can 
    pull back the normal conformal Cartan connection to the orthonormal frame bundle which will yield a 
    formula in terms of the metric $g$, i.e., in terms of the 
    Levi-Civita connection and Schouten tensor. More precisely, the metric $g$ induces a 
    reduction $\iota:\mathcal{P}^g\to\mathcal{P}^0$, and the Levi-Civita connection form
    $A^g\in\Omega^1(\mathcal{P}^g,\mathfrak{so}(p,q))$ determines a $B_0-$equivariant 
    section $\sigma^g:\mathcal{P}^0\to \mathcal{P}^1$ by 
    $\sigma^g(u):=\ker(\gamma^g_u)$, where $\gamma^g$ is the 
    extension of $A^g$ to the conformal frame bundle. Then we have 
    \begin{align*}
      (\sigma^g\circ\iota)^*w^{nc}_s(Y)=[s]^{-1}(d\pi^g_s(Y))+A^g_s(Y)
       -\sum_{i=1}^nP^g_{\pi^g(s)}(d\pi^g_u(Y),s_i)\cdot e_i^*,
    \end{align*}
    where $\pi^g:\mathcal{P}^g\to M$ is the projection map, 
    $s\in\mathcal{P}^g$, $[s]:\mathfrak{b}_{-1}\to T_{\pi^g(s)}M$ 
    the induced isomorphism from $TM\simeq \mathcal{P}^g\times_{(O(p,q),Ad)}\mathfrak{b}_{-1}$, 
    $Y\in T_s\mathcal{P}^g$, $P^g$ 
    denotes the Schouten tensor with respect to $g$, 
    and $\{e_i\}$ is an orthonormal basis in $\mathfrak{b}_{-1}\simeq\R^{p,q}$ 
    with dual basis $\{e_i^*\}$ in $\mathfrak{b}_1$ such that 
    $\{s_i:=[s,e_i]\}$ is an orthonormal basis for $T_{\pi^g(s)}M$ with respect to $g$.

    Now we will define the first prolongation of a conformal spin structure. This requires the pull back, denoted by 
    $\tilde{\cdot}$, of the groups $G, B, B_0$ and $B_1$ by the covering 
    map $\lambda:Spin(p+1,q+1)\to SO(p+1,q+1)$. Consider a conformal spin structure $(\mathcal{Q}^0,f^0)$ 
    on $(M,c)$ and define the set
    \begin{align*}
     \mathcal{Q}^1:=\{\tilde{H}_q\subset T_q\mathcal{Q}^0\mid q\in\mathcal{Q}^0, df^0_q(\tilde{H}_q)\in\mathcal{P}^1\},
    \end{align*}
    and a $\tilde{B}-$action on it by
    \begin{align*}
     \tilde{H}_q\cdot\tilde{b}:=(df^0_{q\cdot b_0})^{-1}\left(df^0_q(\tilde{H}_q)\cdot\lambda(\tilde{b})\right),
    \end{align*}
    for $\tilde{H}_q\in\mathcal{Q}^1$, and $\tilde{b}=\tilde{b}_0\cdot\tilde{b}_1\in\tilde{B}$ 
    ($\tilde{B}$ inherits the semi dirct product structure from $B$). With the 
    obvious projection maps this gives us a $\tilde{B}_1-$principal bundle 
    $(\mathcal{Q}^1,\tilde{\pi}^1,\mathcal{Q}^0,\tilde{B}_1)$, and a $\tilde{B}-$principal bundle 
    $(\mathcal{Q}^1,\tilde{\pi}^0,M,\tilde{B})$ equipped with an equivariant bundle map 
    $f^1:=df^0:\mathcal{Q}^1\to\mathcal{P}^1$. Hence, $(\mathcal{Q}^1,f^1)$ is referred to as the 
    first prolongation of the conformal spin structure. We can lift the normal conformal Cartan connection 
    $w^{nc}$ to a Cartan connection 
    $\tilde{w}^{nc}:=\lambda_*\circ w^{nc}\circ df^1\in\Omega^1(\mathcal{Q}^1,\mathfrak{spin}(p+1,q+1))$ 
    on $\mathcal{Q}^1$. Again, a choice of a metric $g$ from the conformal class leads to the spin connection form 
    $\tilde{A}^g\in\Omega^1(\mathcal{Q}^g,\mathfrak{spin}(p,q))$ which extends to a connection form 
    $\tilde{\gamma}^g$ on $\mathcal{Q}^0$. This in turn induces a $\tilde{B}_0-$equivariant section 
    $\tilde{\sigma}^g:\mathcal{Q}^0\to\mathcal{Q}^1$. Using the reduction 
    map $\tilde{\iota}:\mathcal{Q}^g\to\mathcal{Q}^0$ 
    the pull back of $\tilde{w}^{nc}$ by $\tilde{\sigma}^g\circ\tilde{\iota}$ gives us
    \begin{align*}
     (\tilde{\sigma}^g\circ\tilde{\iota})^*\tilde{w}^{nc}_q(\tilde{Y})
      =\lambda_*^{-1}\bigg(&[f^g(q)]^{-1}d\pi^g_{f^g(q)}(Y)+A^g_{f^g(q)}(Y)\\
     &-\sum_i P^g_{\pi^g\circ f^g(q)}(d\pi^g_{f^g(q)}(Y),s_i)\cdot e_i^*\bigg),
    \end{align*}
    where $\pi^g:\mathcal{P}^g\to M$ is the projection, $\tilde{Y}\in T_q\mathcal{Q}^g$, $Y:=df^g_q(\tilde{Y})$, 
    $P^g$ denotes the Schouten tensor with respect to $g$, and $\{e_i\}$ and $\{e_i^*\}$ are as above such that 
    $\{s_i:=[f^g(q),e_i]\}$ is an orthonormal basis for $T_{\pi^g\circ f^g(q)}M$ with respect to $g$. 

    Summarizing, we have defined first prolongations for the conformal frame bundle and the conformal spin 
    structure of $(M,c)$, and equipped them with distinguished Cartan connections. 
    These structures are the analogues of the orthonormal 
    frame bundle equipped with the Levi-Civita connection fomr and the spin connection form, 
    for a chosen spin structure. 

  \subsection{Tractor bundles for conformal spin structures}\label{tractorbundles}
    Let $(M,c)$ be a conformal spin manifold and $\mathcal{P}^1$ and $\mathcal{Q}^1$ their 
    associated $B-$ and $\tilde{B}-$principal bundles. Considering the standard representation
    $\rho:SO(p+1,q+1)\to Gl(n+2,\R)$ and spin representation 
    $\tilde{\rho}:=\kappa_{p+1,q+1}:Spin(p+1,q+1)\to Gl(\Delta_{p+1,q+1})$, 
    we may define the standard tractor bundle and spin tractor bundle by
    \begin{align*}
      \mathcal{T}(M)&:=\mathcal{P}^1\times_{(B,\rho)}\R^{n+2},\\
      \mathcal{S}(M)&:=\mathcal{Q}^1\times_{(\tilde{B}_0,\tilde{\rho})}\Delta_{p+1,q+1},
    \end{align*}
    where the subscript $\cdot_0$ denotes the connected component of $\tilde{B}$ containing the identity. 
    Both bundles can be equipped with a bundle metric, defined by 
    $g^{\mathcal{T}}(t_1,t_2):=\langle y_1,y_2\rangle_{p+1,q+1}$, 
    for $t_i=[H,y_i]\in\mathcal{T}(M)$, $i=1,2$; and 
    $g^{\mathcal{S}}(s_1,s_2):=(v_1,v_2)_\Delta$, 
    for $s_i=[\tilde{H},v_i]\in\mathcal{S}(M)$, $i=1,2$, since 
    $\langle\cdot,\cdot\rangle_{p+1,q+1}$ and $(\cdot,\cdot)_\Delta$ are invariant under $B$ and $\tilde{B}_0$. 
    Since we have used representations of the groups $SO(p+1,q+1)$ and $Spin(p+1,q+1)$ to form the associated 
    vector bundles, we may define covariant derivatives $\nabla^\mathcal{T}$ and $\nabla^\mathcal{S}$ 
    induced by the Cartan connections $w^{nc}$ and $\tilde{w}^{nc}$. It turns 
    out that $g^\mathcal{T}$ and $g^\mathcal{S}$ 
    are parallel with respect to the corresponding covariant derivatives. 

    Choosing a metric $g$ from the conformal class, the orthonormal frame bundle $\mathcal{P}^g$ is a 
    $SO(p,q)\to CO(p,q)-$reduction of the conformal frame bundle 
    $\mathcal{P}^0$, and a $SO(p,q)\to B-$reduction of 
    the first prolongation $\mathcal{P}^1$. Similarly, 
    $\mathcal{Q}^g$ is a $Spin(p,q)\to CSpin(p,q)-$reduction of the 
    conformal spin structure $\mathcal{Q}^0$, and a $Spin(p,q)\to\tilde{B}-$reduction 
    of the first prolongation $\mathcal{Q}^1$. 
    Thus the following isomorphisms arise:
    \begin{align*}
      \mathcal{T}(M)\simeq& \mathcal{P}^g\times_{(O(p,q),\rho)}\R^{n+2},\\
      \mathcal{S}(M)\simeq&  \mathcal{Q}^g\times_{(Spin_0(p,q),\tilde{\rho})}\Delta_{p+1,q+1},\\
      TM\simeq& \mathcal{P}^g\times_{(O(p,q),Ad)}\mathfrak{b}_{-1}
        \simeq \mathcal{Q}^g\times_{(Spin_0(p,q),Ad\circ\lambda)}\mathfrak{b}_{-1},\\
      T^*M\simeq& \mathcal{P}^g\times_{(O(p,q),Ad)}\mathfrak{b}_1
        \simeq\mathcal{Q}^g\times_{(Spin_0(p,q),Ad\circ\lambda)}\mathfrak{b}_1,\\
      \mathfrak{so}(TM,g)\simeq& \mathcal{P}^g\times_{(O(p,q),Ad)}\mathfrak{so}(p,q)
        \simeq\mathcal{Q}^g\times_{(Spin_0(p,q),Ad\circ\lambda)}\mathfrak{so}(p,q).
    \end{align*}
    Therefore, for $V$ being one of the bundles $TM,T^*M$ or $\mathfrak{so}(TM,g)$, 
    we may define actions $\rho^g:V\to End(\mathcal{T}(M))$ and $\tilde{\rho}^g:V\to End(\mathcal{S}(M))$ 
    by
    \begin{align*}
      \rho^g(\Theta)t:=&[e,\rho_*([e]^{-1}\Theta)y],\\
      \tilde{\rho}^g(\Theta)s:=&[q,\tilde{\rho}_*\circ\lambda_*^{-1}([q]^{-1}\Theta)v],
    \end{align*}
    where $t=[e,y]\in\mathcal{T}(M)$, $s=[q,v]\in\mathcal{S}(M)$, $\Theta\in V$, and $[e]:W\to V$ 
    and $[q]:W\to V$ are the induced isomorphisms $[e]w:=[e,w]$ and $[q]w:=[q,w]$, 
    for $w\in W=\mathfrak{b}_{-1},\mathfrak{b}_1,\mathfrak{so}(p,q)$, respectively. 
    In terms of these actions we have
    \begin{align*}
      \nabla^\mathcal{T}_Xt=&\nabla^g_Xt+\rho^g(X)t-\rho^g(P^g(X))t,\\
      \nabla^\mathcal{S}_Xs=&\nabla^g_Xs+\tilde{\rho}^g(X)s-\tilde{\rho}^g(P^g(X))s,
    \end{align*}
    for sections $t=[e,y]\in\Gamma(\mathcal{T}(M))$, $s=[q,v]\in\Gamma(\mathcal{S}(M))$, and a 
    vector field $X\in\mathfrak{X}(M)$. 
    Note that $\nabla^g_Xt$ and $\nabla^g_Xs$ are abbreviations for $[e,X(y)+\rho_*(A^g_e(de(X)))y]$ 
    and $[q,X(v)+\tilde{\rho}_*(\tilde{A}^g_q(dq(X)))v]$, and $P^g(X)$ is 
    considered as a $1-$form.  

    A crucial step in this subsection is to define a $g-$metric decomposition of standard 
    tractors and spin tractors with respect to a 
    metric $g$ from the conformal class. 
    Firstly, we have the bundle isomorphism
    \begin{align}
      \Phi^g:\mathcal{T}(M)&\to \underline{M}\oplus TM\oplus \underline{M}=:\mathcal{T}(M)_g,\notag\\
      t=[e,y]&\mapsto (\alpha,X,\beta)=:t_g,\label{eq:tractorisom}
    \end{align}
    where $\underline{M}:=M\times\R$ is the trivial bundle, $y\in\R^{n+2}$ has 
    coordinates $(\alpha,x=(x_1,\ldots,x_n),\beta)$ 
    with respect to the basis $\{f_-,e_i,f_+\}$ of $\R^{n+2}$, and $X:=[e]^{-1}x\in TM$.
    Secondly, we have the bundle isomorphism 
    \begin{align}
      \Psi^g:\mathcal{S}(M)&\to S(M,g)\oplus S(M,g)=:\mathcal{S}(M)_g\notag\\
      s=[q,v]&\mapsto (\psi,\phi)=:s_g,\label{eq:spintractorisom}
    \end{align}
    where $\psi=[q,w_1]$ and $\phi=[q,w_2]$, with $w_1,w_2\in\Delta_{p,q}$ being 
    determined as follows: Consider the two 
    $Spin(p,q)-$invariant subspaces $W^\pm:=\{  v\in\Delta_{p+1,q+1}\mid f_\pm\cdot v=0\}$ 
    of $\Delta_{p+1,q+1}$. 
    Note that we naturally identify $W^+$ with $\Delta_{p,q}$. 
    Hence, $\tilde{\rho}$ restricted to $Spin(p,q)$ decomposes into two representations 
    $\tilde{\rho}^\pm:Spin(p,q)\to Gl(W^\pm)$, such 
    that $\tilde{\rho}_{\mid Spin(p,q)}=\tilde{\rho}^+\oplus\tilde{\rho}^-$. 
    From the definition of $W^\pm$ it follows that $\tilde{\rho}^\pm$ are equivalent with respect to the isomorphism 
    $W^+\ni w\mapsto f_-\cdot w\in W^-$. 
    Therefore, our element in question $v\in \Delta_{p+1,q+1}$ can be uniquely decomposed 
    as $v=w_1+f_-\cdot w_2$ with 
    $w_1,w_2\in W^+$, due to the isomorphism 
    $W^+\times W^+\ni(w_1,w_2)\mapsto w_1+f_- \cdot w_2\in\Delta_{p+1,q+1}$. 
    
    With the help of the two maps $\Phi^g$ and $\Psi^g$ we will interpret tractor objects with 
    data coming from the metric $g$. For example, we have that
    \begin{align*}
      \Phi^g\circ\nabla^{\mathcal{T}}_X\circ(\Phi^g)^{-1}
        =\begin{pmatrix}\nabla^{LC}_X& -P^g(X,\cdot)& 0\\ X\cdot & \nabla^{LC}_X& P^g(X)^\natural\cdot\\ 
             0& -g(X,\cdot)& \nabla^{LC}_X
           \end{pmatrix}
    \end{align*}
    and 
    \begin{align*}
      \Psi^g\circ\nabla^{\mathcal{S}}_X\circ(\Psi^g)^{-1}
        =\begin{pmatrix}\nabla^{S(M,g)}_X& X\cdot\\ \frac 12P^g(X)^\natural\cdot& \nabla^{S(M,g)}_X
          \end{pmatrix},
    \end{align*}
    which follows from the actions $\rho^g$ and $\tilde{\rho}^g$ defined above. 
    A further example is given by the bundle metrics, here we have that
    \begin{align}
      g^\mathcal{T}(t_1,t_2)=\alpha_1\beta_2+g(X_1,X_2)+\beta_2\alpha_1,\label{eq:tractormetric}
    \end{align}
    for $t_i=[e,y_i]$, $i=1,2$. Moreover, for $s_i=[q,v_i]\in\mathcal{S}(M)$, $i=1,2$, we have that
    \begin{align}
      g^\mathcal{S}(s_1,s_2)=-2\sqrt{2}i^p\left(<\phi_1,\psi_2>
        +(-1)^p<\psi_1,\phi_2>\right).\label{eq:spintractormetric}
    \end{align}
    Note that these results are based on the isomorphisms \eqref{eq:tractorisom} and \eqref{eq:spintractorisom}.
    Let us end this subsection with the realization of standard and spin tractors with respect to two metrics 
    $g$ and $\hat{g}=e^{2\sigma}g$ from the conformal class. Here it holds that
    \begin{align*}
      T(g,\sigma):=\Phi^{\hat{g}}\circ(\Phi^g)^{-1}=
        \begin{pmatrix}e^{-\sigma}& -e^{-\sigma}d\sigma& -\frac 12 e^{-\sigma}\abs{\grad^g(\sigma)}_g^2\\
          0&e^{-\sigma}&e^{-\sigma}\grad^g(\sigma)\\ 0&0&e^{\sigma}  \end{pmatrix}
    \end{align*}
    and 
    \begin{align*}
      T^{\mathcal{S}(M)}(g,\sigma):=\Psi^{\hat{g}}\circ(\Psi^g)^{-1}=F_\sigma\oplus F_\sigma
        \begin{pmatrix}e^{\frac 12\sigma}& 0\\
          \frac 12e^{-\frac 12\sigma}\grad^g(\sigma)\cdot&e^{-\frac 12\sigma}  \end{pmatrix},
    \end{align*}
    where $F_\sigma:S(M,g)\to S(M,\hat{g})$ is the bundle isomorphism relating spinor bundles 
    for two conformally related metrics $g$ and $\hat{g}=e^{2\sigma}g$.

\section{Relevant differential operators}\label{relevantoperators}
   
  In this section we present some operators necessarily for the construction of 
  conformal powers of the Dirac operator. First we recall the 
  construction of the splitting operator 
  for the standard tractor bundle (in the spirit of \cite{CapSoucek, CapGoverSoucek}), 
  and compute its formal adjoint. The notation is 
  borrowed from these two papers. Both the splitting operator and its adjoint can be 
  extended to $\mathcal{S}(M)$, as well as to 
  $\mathcal{S}^k(M):=\otimes^k(\mathcal{T}(M))\otimes\mathcal{S}(M)$, for $k\geq 0$. Secondly, we 
  consider the translations of the strongly invariant Yamabe operator with these splitting operators and 
  their formal adjoints. We do this in order 
  to obtain higher order differential operators acting on the spin tractor bundle.  

  Let us assume that $M$ is even dimensional, and so, that $p+1+q+1=2(m+1)$. 
  The odd dimensional case is treated similarly. 
  The weighted standard tractor bundle $\mathcal{T}(M)[w-1]$ splits under the conformal group $B_0$ as 
  $\mathcal{T}(M)_g[w-1]=\underline{M}[w-2]\oplus TM[w-1]\oplus \underline{M}[w]$. The lowest 
  weights for these summands are $(w-2|0,\ldots,0)$, $(w-1|1,0,\ldots,0)$ and $(w|0,\ldots,0)$, each of 
  length $(m+1)$. Moreover, we denote by 
  $\rho=(m,m-1,\ldots,1,0)$ the half sum of all positive roots. 

  The curved Casimir operator $C:\Gamma(\mathcal{T}(M)_g)\to \Gamma(\mathcal{T}(M)_g)$ 
  obeys the following formula, given in \cite[Section $2.2$]{CapGoverSoucek},
  \begin{align}
      C(t_g)=\beta(t_g)-2\sum_{l=1}^n\rho^g(\xi_l)
       \left(\nabla^g_{\xi^l}t_g-\rho^g(P(\xi^l))t_g\right)\label{eq:Casimir},
  \end{align}
  where $\{\xi^l\}$ denotes a basis of $TM$ and $\{\xi_l\}$ is its 
  dual, $t_g\in\Gamma(\mathcal{T}(M)_g)$, $P(\xi^l)$ 
  is considered to be a $1-$form, and the 
  map $\beta:\Gamma(\mathcal{T}(M)_g)\to\Gamma(\mathcal{T}(M)_g)$ 
  acts on the direct sum by the Casimir scalars 
  \begin{align*}
      \beta_1=w(w+n)-2(n+2w-2),\quad \beta_2=w(w+n)-2w,\quad \beta_3=w(w+n),
  \end{align*}
  which can be derived from \cite[Theorem $1$]{CapSoucek}. Thus, using 
  \begin{align*}
      C(t_g)=\begin{pmatrix}
                 [w(w+n)-2(n+2w-2)]\alpha-2\diver(X)-2J\beta\\ [w(w+n)-2w]X+2\grad(\beta)\\ w(w+n)\beta
              \end{pmatrix},
  \end{align*}
  where $t_g=(\alpha,X,\beta)\in\Gamma(\mathcal{T}(M)_g)$, one computes that
  \begin{align*}
       (C-\beta_1)\circ(C-\beta_2)t_g=4\begin{pmatrix}
                                                     -\Delta_g^{\nabla^{LC}}\beta-wJ\beta\\ 
                                                     (n+2w-2)(d\beta)^\natural\\ 
                                                     w(n+2w-2)\beta 
                                                   \end{pmatrix}=:4D(g)\beta.
  \end{align*}
  Note our sign convention for the Laplacian is $\Delta_g^\nabla:=Tr_g(\nabla^{LC}\circ\nabla^{LC})$. 
  This defines a mapping $D(g):\Gamma(\underline{M}[w])\to \Gamma(\mathcal{T}(M)_g[w-1])$. 
  In the same manner one 
  constructs an operator 
  $D^{k}:\Gamma(\mathcal{S}^k(M))\to\Gamma(\mathcal{S}^k(M)\otimes\mathcal{T}(M)_g)$, 
  for $k\geq 0$. The splitting operator for the spinor bundle is constructed similarly: 
  The spin tractor bundle splits under the conformal spin group $\tilde{B}_0$ as 
  $\mathcal{S}(M)_g\simeq S(M,g)[\frac 12]\oplus S(M,g)[-\frac 12]$. Thus, $\mathcal{S}(M)_g[\eta-\frac 12]$ 
  decomposes 
  into a direct sum corresponding to lowest weights $(\eta|\frac 12,\ldots,\frac 12)$ 
  and $(\eta-1|\frac 12,\ldots,-\frac 12)$. Again, the Casimir scalars are given by 
  \begin{align*}
      \beta_1=\eta(\eta+n)+\frac 12\sum_{i=0}^{m-1}(\frac 12 +2i),\quad 
        \beta_2=(\eta-1)(\eta-1+n)+\frac 12\sum_{i=0}^{m-1}(\frac 12+2i).
  \end{align*}
  Hence, equation \eqref{eq:Casimir} adapted to the spin tractor setting gives us 
  \begin{align*}
       C(s)=\begin{pmatrix}  \beta_1\psi\\ \beta_2\phi+\slashed{D}\psi\end{pmatrix},
  \end{align*}
  for $s=(\psi,\phi)\in\Gamma(\mathcal{S}(M)_g)$. This shows that
  \begin{align*}
      (C-\beta_2)s=2\begin{pmatrix} (\eta+\frac{n-1}{2})\psi\\ 
                              \frac 12\slashed{D}\psi 
                            \end{pmatrix}=:2D^{spin}(g)\psi
  \end{align*}
  defines a map $D^{spin}(g):\Gamma(S(M,g)[\eta])\to\Gamma(\mathcal{S}(M)_g[\eta-\frac 12])$. Note that the
  construction of $D^k(g)$ and $D^{spin}(g)$ only depends on the tractor data, hence they are well defined.

  From now on we will work with unweighted bundles. The conformal weights are absorbed into the 
  splitting operators as follows: 
  \begin{align}
    D^{k}(g,w):\Gamma(\mathcal{S}^k(M))&\to \Gamma(\mathcal{S}^k(M)\otimes\mathcal{T}(M)_g)\notag\\
    s&\mapsto \begin{pmatrix} -\Box_w^\nabla s \\
                        (n-2+2w)(\nabla s)^\natural\\
                        w(n-2+2w)s 
                    \end{pmatrix},\label{eq:TractorD}
  \end{align}
  where $\Box_w^\nabla s:=\Delta_g^\nabla s+wJs$, and 
  \begin{align*}
    D^{spin}(g,\eta):\Gamma(S(M,g))&\to \Gamma(\mathcal{S}(M)_g)\notag\\
    \psi&\mapsto \begin{pmatrix} (\eta+\frac{n-1}{2})\psi\\ \frac 12\slashed{D}\psi
                        \end{pmatrix}.
  \end{align*}
  Since we are restricting our attention to unweighted bundles we have the 
  following conformal transformation laws:
  \begin{sat}\label{TractorD}
    Let $\hat{g}=e^{2\sigma}g$, $s\in\Gamma(\mathcal{S}^k(M))$ and $\psi\in\Gamma(S(M,g))$. Then one has
    \begin{align*}
      D^k(\hat{g},w)(e^{w\sigma}s)=&e^{(w-1)\sigma} T(g,\sigma)D^k(g,w)s,\\
      D^{spin}(\hat{g},\eta)(e^{\eta\sigma}F_\sigma\psi)
        =&e^{(\eta-\frac 12)\sigma}T^{\mathcal{S}(M)}(g,\sigma)D^{spin}(g,\eta)\psi,
    \end{align*}
    for all $w,\eta\in\R$. Here $F_\sigma:\Gamma(S(M,g))\to\Gamma(S(M,\hat{g}))$ is the 
    isomorphism for conformally related metrics. 
  \end{sat}
  For later purposes, let us define
  \begin{align*}
    C^{spin}(g,\eta):\Gamma(\mathcal{S}(M)_g)&\to \Gamma(S(M,g))\\
    s_g=(\psi,\phi)&\mapsto\frac 12\slashed{D}\psi-(\eta+\frac n2)\phi,
  \end{align*}
  and 
  \begin{align}
    C^{k}(g,w): \Gamma(\mathcal{S}^k(M)\otimes\mathcal{T}(M)_g)\to& \Gamma(\mathcal{S}^k(M))\notag\\
    (s_1,\eta,s_2)\mapsto& (n+nw_1+w_1w)s_1+(n+2w)\diver(\eta)\notag\\
    &-(\Delta_g^\nabla+(1-n-w)J)s_2,\label{eq:TractorC}
  \end{align}
  where $\diver(Y\otimes s):=\diver(Y) s+\nabla_Ys\in \Gamma(\mathcal{S}^k(M))$, for 
  $Y\otimes s\in\Gamma(TM\otimes\mathcal{S}^k(M))$, and the divergence of a vector field is defined by 
  \begin{align*}
    \diver(Y):=\sum_i\varepsilon_ig(\nabla_{s_i}Y,s_i),
  \end{align*}
  in terms of a local section 
  $(s_1,\ldots s_n):U\subset M\to\mathcal{P}^g$. By the proposition below, 
  they are the formal adjoints of corresponding splitting operators. 
  \begin{sat}\label{SelfadjointTractorD}
    As formal adjoints with respect to the corresponding $L^2-$scalar product we have that 
    \begin{align*}
      \left(D^k(g,w)\right)^*&=C^k(g,1-n-w),\\
      \left(C^k(g,w)\right)^*&=D^k(g,1-n-w),\\
      \left(D^{spin}(g,\eta)\right)^*&=-2\sqrt{2}i^pC^{spin}(g,\frac 12-n-\eta),\\
      \left(C^{spin}(g,\eta)\right)^*&=-\frac{1}{2\sqrt{2}}i^pD^{spin}(g,\frac 12-n-\eta).
    \end{align*} 
  \end{sat}
  \begin{bew}
    Using the formulas \eqref{eq:tractormetric} 
    and \eqref{eq:spintractormetric} for the scalar products $g^\mathcal{T}$ and $g^\mathcal{S}$ 
    we compute, for $k=0$, that 
    \begin{align*}
       g^{\mathcal{T}\otimes\mathcal{S}}&\left(D^0(g,w)s,\begin{pmatrix} s_1\\ 
                \eta\\ s_2\end{pmatrix}\right)_{L^2}\\
         =&-g^\mathcal{S}(\Box_w^\nabla s,s_2)_{L^2}
         +g^{TM\otimes \mathcal{S}}\left(w_1(\nabla s)^\natural,\eta\right)_{L^2}
         +g^\mathcal{S}(ww_1s,s_1)_{L^2}\\
       =&-g^\mathcal{S}(s,\Box_w^\nabla s_2)_{L^2}
         -g^{TM\otimes \mathcal{S}}\left(s,w_1\diver(\eta)\right)_{L^2}
         +g^\mathcal{S}(s,ww_1s_1)_{L^2}\\
       =&g^\mathcal{S}\left(s, C^0(g,1-n-w)\begin{pmatrix}s_1\\ \eta\\ s_2  \end{pmatrix}\right)_{L^2},
    \end{align*}
    where we have used the known adjoints of $\Delta_g^\nabla$ and $d^\nabla$. 
    Note that the index $\cdot_{L^2}$ indicates the induced $L^2-$scalar product. 
    The case for $k>0$ runs along the same lines. 
    The second assertion follows immediately. Coming to the third one, 
    we have, for $\eta_1:=(\eta+\frac{n-1}{2})$, that 
    \begin{align*}
      g^\mathcal{S}&\left(D^{spin}(g,\eta)\psi,\begin{pmatrix} \phi_1\\ \phi_2\end{pmatrix}\right)_{L^2}\\
      =&-2\sqrt{2}i^p\left( <\frac 12\slashed{D}\psi,\phi_1>_{L^2}+(-1)^p<\eta_1\psi,\phi_2>_{L^2}\right)\\
      =&-2\sqrt{2}i^p\left( (-1)^p<\psi,\frac 12\slashed{D}\phi_1>_{L^2}
        +(-1)^p<\psi,\eta_1\phi_2>_{L^2}\right)\\
      =&-2\sqrt{2}i^p(-1)^p\left(<\psi,\frac 12\slashed{D}\phi_1
        -(\frac 12-\eta-n+\frac{n}{2})\phi_2>_{L^2}\right)\\
      =&<\psi, (-2)\sqrt{2}i^pC^{spin}(g,\frac 12-n-\eta)\begin{pmatrix} \phi_1\\ \phi_2\end{pmatrix}>_{L^2},
    \end{align*}
    where we haved used the (anti-) self-adjointness of $\slashed{D}$. 
    Also note the hermiticity of $<\cdot,\cdot>_{L^2}$. An analogous computation shows that
    \begin{align*}
      <C^{spin}(g,\eta)\begin{pmatrix}\phi_1\\ \phi_2 \end{pmatrix},\psi>_{L^2}
        =&<\frac 12\slashed{D}\phi_1,\psi>_{L^2}-<(\eta+\frac n2)\phi_2,\psi>_{L^2}\\
      =&-\frac{1}{2\sqrt{2}i^p}g^\mathcal{S}\left(\begin{pmatrix} \phi_1\\ \phi_2\end{pmatrix},
        \begin{pmatrix}(\frac 12-n-\eta+\frac{n-1}{2})\psi\\ \frac 12\slashed{D}\psi  \end{pmatrix}\right)_{L^2}\\
      =&g^\mathcal{S}\left(\begin{pmatrix} \phi_1\\ \phi_2\end{pmatrix},
        \frac{-1}{2\sqrt{2}}i^pD^{spin}(g,\frac 12-n-\eta)\psi\right)_{L^2},
    \end{align*}
    which completes the proof.     
  \end{bew}

  It now follows from this proposition and from the invariance of the corresponding scalar products with respect 
  to $g$ and $\hat{g}=e^{2\sigma}g$, that: 
  \begin{sat}\label{TractorC}
    For $\hat{g}=e^{2\sigma}g$, $s_g=(\psi,\phi)\in\Gamma(\mathcal{S}(M)_g)$ and 
   $(s_1,\eta,s_2)\in\Gamma(\mathcal{S}^k(M)\otimes\mathcal{T}(M)_g)$, one has
    \begin{align*}
      C^k(\hat{g},w)(e^{w\sigma}T(g,\sigma)(s_1,\eta,s_2))=&e^{(w-1)\sigma} C^k(g,w)(s_1,\eta,s_2),\\
      C^{spin}(\hat{g},\eta)(e^{\eta\sigma}T^{\mathcal{S}(M)}(g,\sigma)s_g)
        =&e^{(\eta-\frac 12)\sigma}F_\sigma(C^{spin}(g,\eta)s_g).
    \end{align*}
  \end{sat}

  Since Proposition \ref{TractorD} holds for any real numbers $w,\eta\in\R$, the conformal covariance of the 
  Box operator $\Box^\nabla_{\frac{2-n}{2}}$ (the strongly invariant Yamabe operator), 
  and the Dirac operator $\slashed{D}$, follow from that of 
  $D^k(g,w)$ and $D^{spin}(g,\eta)$, for the values $w=\frac{2-n}{2}$ and $\eta=\frac{n-1}{2}$. 
 
  As mentioned above, the operator $\Box^\nabla_{\frac{2-n}{2}}$ acts conformally 
  on $\Gamma(\mathcal{S}^k(M))$, for $k\in\N_0$. 
  Hence we can use the curved translation principle, introduced in \cite{EastwoodRice}, to define 
  $P^{\mathcal{S}(M)}_{2N}(g):\Gamma(\mathcal{S}(M)_g)\to\Gamma(\mathcal{S}(M)_g)$ for $N\in\N$, by
  \begin{align}
    P^{\mathcal{S}(M)}_2(g):=&\Box^\nabla_{\frac{2-n}{2}}\notag\\
    P^{\mathcal{S}(M)}_{2N}(g):=&C^0(g,-\frac{2(N-1)+n}{2})\circ\ldots\circ C^{N-2}(g,-\frac{2+n}{2})
     \circ\Box^\nabla_{\frac{2-n}{2}}\circ\notag\\
    &\circ D^{N-2}(g,\frac{4-n}{2})\circ\ldots\circ D^0(g,\frac{2N-n}{2}),\quad N>1.\label{eq:DieP}
  \end{align}
  These operators satisfy the following:
  \begin{sat}\label{ConfCovOfP}
    The operator $P^{\mathcal{S}(M)}_{2N}(g)$ is conformally covariant of 
    bi-degree $(\frac{2N-n}{2},-\frac{2N+n}{2})$, i.e., 
    for $\hat{g}=e^{2\sigma}g$ we have
    \begin{align*}
      P_{2N}^{\mathcal{S}(M)}(\hat{g})(e^{\frac{2N-n}{2}\sigma}s_{\hat{g}})
       =e^{-\frac{2N+n}{2}\sigma}P_{2N}^{\mathcal{S}(M)}(g)s_g,
    \end{align*}
    for $s\in\Gamma(\mathcal{S}(M))$. Its leading term is given 
    by $c(n,N)(\Delta_g^\nabla)^N$, where the constant is 
    \begin{align}
      c(n,N):=(-1)^{N-1}\prod_{k=1}^{N-1}[k(2+2k-n)].\label{eq:constant}
    \end{align} 
  \end{sat}
  \begin{bew}
    The conformal covariance follows from the well-chosenness of $w$ in the 
    composition. The given expression for $c(n,N)$ follows directly from \ref{eq:TractorD}, producing $(-1)^{N-1}$, 
    and \ref{eq:TractorC}, producing the product.   
  \end{bew}

  \begin{bem}
    In case of even $n$, the operator $P^{\mathcal{S}(M)}_{2N}(g)$, for $N\geq \frac n2$, 
    is not identically zero as stated in \cite[Proposition $5.26$]{Fischmann}. It is just of order less than $2N$, 
    due to the fact that the constant $c(n,N)$ is zero in this case. 
  \end{bem}

\section{The construction of conformal powers of the Dirac operator and related structures}\label{results}
  This section makes further use of the curved translation principle, \cite{EastwoodRice}, to define  
  conformally covariant operators, acting on the spinor bundle, which are conformal powers of the Dirac operator. 
  Furthermore, we present explicit formulas for lower order examples in general, and 
  subsequently simplify to the Einstein case. We then go on to prove some formal self-adjointness 
  results. Using these explicit formulas we are able to show that 
  the conformal powers of the Dirac operator are polynomials in first order operators. 

  Consider the differential operator   
  \begin{align}
    D_{2N+1}(g):=C^{spin}(g,-\frac{2N+n}{2})\circ P^{\mathcal{S}(M)}_{2N}(g)
      \circ D^{spin}(g,\frac{2N+1-n}{2})\label{eq:ConfDirac}
  \end{align}
  constructed from $P^{\mathcal{S}(M)}_{2N}(g)$ by translation, which acts on the spinor bundle. 
  \begin{theorem}
    Let $N\in\N$. The operator $D_{2N+1}(g)$ is conformally covariant of bi-degree 
    $(\frac{2N+1-n}{2},-\frac{2N+1+n}{2})$, i.e., 
    for $\hat{g}=e^{2\sigma}g$ and $\psi\in\Gamma(S(M,g))$ we have 
    \begin{align*}
       D_{2N+1}(\hat{g})(e^{\frac{2N+1-n}{2}\sigma} F_\sigma\psi)
        =e^{-\frac{2N+1+n}{2}\sigma}F_\sigma\circ D_{2N+1}(g)\psi.
    \end{align*}
    Its leading term is given by a constant multiple of 
    $\slashed{D}^{2N+1}$.
  \end{theorem}
  \begin{bew}
    The conformal covariance follows directly from the construction of $D_{2N+1}(g)$. 
    The leading term is given by a scalar multiple of $\slashed{D}^{2N+1}$, due to the fact that 
    $P^{\mathcal{S}(M)}_{2N}(g)$ has leading term $c(n,N)(\Delta_g^\nabla)^N$ and the explicit formula 
    \begin{align*}
      \Delta_g^\nabla=\begin{pmatrix}-\slashed{D}^2+\frac{n-2}{2}J& 2\slashed{D}\\
                                                      (P,\nabla^{S(M,g)})+\frac 12\grad(J)&-\slashed{D}^2+\frac{n-2}{2}J
                                \end{pmatrix}.
    \end{align*} 
    The scalar multiple  
    of $\slashed{D}^{2N+1}$ is a product of $c(n,N)$ and a term independently of $n$. 
  \end{bew}

  \begin{bem}
    In case of even $n$, the operator $D_{2N+1}(g)$, for $N\geq \frac n2$, 
    is not identically zero as stated in \cite[Theorem $5.27$]{Fischmann}. It is just of order less than $2N+1$, 
    due to the fact that the constant infront of $\slashed{D}^{2N+1}$ is zero in this case. Thus, in that case, 
    the last theorem does not yield conformal powers of the Dirac operator.  
  \end{bem}
  Explicit formulas for $D_{2N+1}(g)$, for $N=1,2$, can be derived from explicit knowledge of 
  $P_{2}^{\mathcal{S}(M)}(g)$ and $P_{4}^{\mathcal{S}(M)}(g)$ found in 
  \cite[Proposition $5.28$]{Fischmann} and \cite[Proposition $5.36$]{Fischmann}:
  \begin{theorem}\label{explicit}
    Let $(M,g)$ be a semi Riemannian spin manifold. The operators $D_3(g)$ and $D_5(g)$ are given by
    \begin{align*}
      D_3(g)=&-\frac 12[\slashed{D}^3-(P,\nabla^{S(M,g)})-(\nabla^{S(M,g)},P\cdot)],\\
      D_5(g)=&(n-4)\Big[\slashed{D}D_3(g)\slashed{D}
        +2(\slashed{D}^2 D_3(g)+D_3(g)\slashed{D}^2)-4\slashed{D}^5\\
      &+4(2P^2+\frac{1}{n-4}B,\nabla^{S(M,g)})+4(\nabla^{S(M,g)},2P^2\cdot+\frac{1}{n-4}B\cdot)\\
      &-2(C,P\cdot)-2(P,C\cdot)\Big]\\
      &+\slashed{D}(W\cdot W\cdot )+W\cdot W\cdot\slashed{D}+4(C\cdot W\cdot+W\cdot C\cdot),
    \end{align*}
    where the bracket and product notations were introduced in Subsection \ref{Preliminaries}.
  \end{theorem}
  A detailed presentation of the proof can be found in 
  \cite[Theorem $5.29$]{Fischmann} and \cite[Theorem $5.39$, Remark $5.40$]{Fischmann}.
  
  We have to remark that the operator $P^{\mathcal{S}(M)}_4(g)$ decomoses 
  into $P_{4}^{\mathcal{S}(M)}(g)=P_4(g)+R(g)$, where 
  both operators are conformally covariant of the same bi-degree as 
  $P_{4}^{\mathcal{S}(M)}(g)$. However, $P_4(g)$ has leading 
  term a multiple of $(\Delta_g^\nabla)^2$, whereas $R(g)$ is a zero order 
  operator involving Weyl and Cotton curvatures, see 
  \cite[Proposition $5.37$, Remark $5.38$]{Fischmann}. Hence, the 
  operator $D_5(g)$ decomposes into $D_5(g)=D^{red}_5(g)+R^{spin}(g)$, where 
  \begin{align*}
    R^{spin}(g):=&C^{spin}(g,-\frac{4+n}{2})\circ R(g)\circ D^{spin}(g,\frac{5-n}{2}) \\
    =&\slashed{D}(W\cdot W\cdot )+W\cdot W\cdot\slashed{D}+4(C\cdot W\cdot+W\cdot C\cdot). 
  \end{align*}
  Terms involving the Weyl curvature in the formula for $D_{5}(g)$ are relics from the tractor 
  machinery we used for the construction. 
  
  Finally, let us denote the first three examples of conformal powers of the Dirac operator as follows:
  \begin{align*}
   \mathcal{D}_1:=\slashed{D};\quad \mathcal{D}_3:=-2D_3(g);
     \quad \mathcal{D}_5:=\frac{1}{n-4}D^{red}_5(g),\quad (n\neq 4).
  \end{align*}
  These operators have an odd power of the Dirac operator as the leading term. 
  Due to the explicit formulas we can prove the following theorem. 
  \begin{theorem}\label{Einstein}
    Let $(M^n,g)$ be an Einstein spin manifold. Then one has
    \begin{align*}
       \mathcal{D}_3=&\left(\slashed{D}-\sqrt{\frac{2J}{n}}\right)\slashed{D}\left(\slashed{D}
         +\sqrt{\frac{2J}{n}}\right),\\
       \mathcal{D}_5=&\left(\slashed{D}-\sqrt{\frac{8J}{n}}\right)\left(\slashed{D}
         -\sqrt{\frac{2J}{n}}\right)\slashed{D}
          \left(\slashed{D}+\sqrt{\frac{2J}{n}}\right)\left(\slashed{D}+\sqrt{\frac{8J}{n}}\right),
    \end{align*}
    where $J$ is the normalized (constant) scalar curvature. 
  \end{theorem}
  \begin{bew}
    Since $(M,g)$ is Einstein, we have by definition that $Ric=\lambda g$, 
    for some constant $\lambda\in\R$. Thus, the scalar curvature 
    satisfies $\tau=n\lambda$, hence it is constant, and so is $J=\frac{n\lambda}{2(n-1)}$. 
    It follows that $P=\frac{\lambda}{2(n-1)}g$. This shows, that
    \begin{align*}
      \mathcal{D}_3=&\slashed{D}^3-2(P,\nabla)=\slashed{D}^3-\frac{\lambda}{n-1}\slashed{D}\\
      =&\left(\slashed{D}-\sqrt{\frac{2J}{n}}\right)\slashed{D}\left(\slashed{D}+\sqrt{\frac{2J}{n}}\right).
    \end{align*}
    Since the Bach tensor and the Cotton tensor vanish for Einstein metrics, we have
    \begin{align*}
      \mathcal{D}_5=&\slashed{D}\mathcal{D}_3\slashed{D}
       +2\left(\slashed{D}^2\mathcal{D}_3+\mathcal{D}_3\slashed{D}^2\right)-4\slashed{D}^5+16(P^2,\nabla)\\
      =&\slashed{D}^5-5\frac{\lambda}{n-1}\slashed{D}^3+4\frac{\lambda^2}{(n-1)^2}\slashed{D}\\
      =&\left(\slashed{D}-\sqrt{\frac{8J}{n}}\right)\left(\slashed{D}-\sqrt{\frac{2J}{n}}\right)\slashed{D}
          \left(\slashed{D}+\sqrt{\frac{2J}{n}}\right)\left(\slashed{D}+\sqrt{\frac{8J}{n}}\right),
    \end{align*}
    which completes the proof.
  \end{bew}

  The result of the last theorem is analogous to the product structure for the conformal 
  powers of the Laplacian for Einstein manifolds, compare \cite{Gover}. 
  For example, Theorem \ref{Einstein} in case of 
  the standard sphere, i.e., $J=\frac n2$,  agrees with the result obtained in \cite{ES}, where it was proven that 
  all conformal odd powers of the Dirac operator have such a product structure.

  In order to prove some formal (anti-) self-adjointness results, we present the following theorem. It 
  generalizes the formal (anti-)
  self-adjointness of the Dirac operator, which is given in terms of the bracket notation \eqref{eq:bracket1} by 
  \begin{align*}
      \slashed{D}=\frac 12\big((g,\nabla^{S(M,g)})+(\nabla^{S(M,g)},g\cdot)\big),
  \end{align*}
  to arbitrary symmetric $(0,2)-$tensor fields $T$ instead of $g$. 
  \begin{theorem}\label{selfadjoint}
      Let $(M,g)$ be a semi Riemannian spin manifold without boundary, and let 
      $T$ be a symmetric $(0,2)-$tensor field. The operator 
      \begin{align*}
        (T,\nabla^{S(M,g)})+(\nabla^{S(M,g)},T\cdot):\Gamma\left(S(M,g)\right)\to\Gamma\left(S(M,g)\right)
      \end{align*}
      is formally self-adjoint, or anti self-adjoint, with respect to the $L^2-$scalar product, depending on the 
      signature $(p,q)$ of $(M,g)$.
  \end{theorem}
  \begin{bew}
      Let $\psi,\phi\in\Gamma_c\left(S(M,g)\right)$ be the compactly supported spinors, and 
      define a $1-$form with values in $\C$ by 
      $w(X):=<T(X)^\natural\cdot\psi,\phi>$. Considering its dual $Y_w$, with respect to $g$, 
      and taking its divergence we obtain
      \begin{align*}
        \diver(Y_w)=&\sum_i\varepsilon_i\big[<T(s_i)^\natural\cdot\nabla^{S(M,g)}_{s_i}\psi,\phi>
          -(-1)^p<\psi,T(s_i)^\natural\cdot\nabla_{s_i}^{S(M,g)}\phi>\big]\\
        &+(-1)^p<\psi,(\delta^{\nabla^{LC}}T)^\natural\cdot\phi>.
      \end{align*}
      Using Stokes' Theorem we get $\int_M\diver(Y_w)dM=0$, hence 
      \begin{align*}
        <(T,\nabla^{S(M,g)}\psi)&+(\nabla^{S(M,g)},T\cdot\psi),\phi>_{L^2}\\
        =&\int_M <(T,\nabla^{S(M,g)}\psi)+(\nabla^{S(M,g)},T\cdot\psi),\phi>dM\\
        =&(-1)^p\int_M<\psi,2(T,\nabla^{S(M,g)}\phi>-(\delta^{\nabla^{LC}}T)^\natural\cdot\phi>dM\\
        =&(-1)^p<\psi,(T,\nabla^{S(M,g)}\phi)+(\nabla^{S(M,g)},T\cdot\phi)>_{L^2},
      \end{align*}
      which completes the proof. 
  \end{bew}

  This leads us to the following result:
  \begin{theorem}
   Let $(M,g)$ be a semi Riemannian spin manifold without boundary. 
   The operators $\mathcal{D}_k$, $k=1,3,5$, are formally self-adjoint 
   (anti self-adjoint) with respect to the $L^2-$scalar product, i.e., 
   \begin{align*}
     <\mathcal{D}_k\psi,\phi>_{L^2}=(-1)^{p}<\psi,\mathcal{D}_k\phi>_{L^2}
   \end{align*}
   for $\psi,\phi$ compactly supported sections of the spinor bundle. 
  \end{theorem}
  \begin{bew}
    This follows from Theorem \ref{selfadjoint}, Theorem \ref{explicit}, and the fact that we have
    \begin{align*}
      <(C,P\cdot)\psi&+(P,C\cdot)\psi,\phi>\\
      =&\sum_i\varepsilon_i<C(s_i)\cdot P(s_i)\cdot\psi+P(s_i)\cdot C(s_i)\cdot\psi,\phi>\\
      =&(-1)^p\sum_i\varepsilon_i<\psi,P(s_i)\cdot C(s_i)\cdot\phi+C(s_i)\cdot P(s_i)\cdot\phi>\\
      =&(-1)^p<\psi, (C,P\cdot)\phi+(P,C\cdot)\phi>,
    \end{align*}
    for any $\psi,\phi\in\Gamma(S(M,g))$, where $\{s_i\}$ is a $g-$orthonormal basis. 
  \end{bew}

  This theorem is a special case of the following result:
  \begin{theorem}
    Let $(M,g)$ be a semi Riemannian spin manifold without boundary. 
    For $N\in\N$ the operator $D_{2N+1}(g)$ is formally self-adjoint (anti self-adjoint) with 
    respect to the $L^2-$scalar product, i.e., 
    \begin{align*}
     <D_{2N+1}(g)\psi,\phi>_{L^2}=(-1)^{p}<\psi,D_{2N+1}(g)\phi>_{L^2}
   \end{align*}
   for $\psi,\phi$ compactly supported sections of the spinor bundle. 
  \end{theorem}
  \begin{bew}
     First of all note that from Proposition \ref{SelfadjointTractorD} the operator 
     $P_{2N}^{\mathcal{S}(M)}(g)$ is formally 
     self-adjoint. Hence, by further use of Proposition \ref{SelfadjointTractorD}, we get that
     \begin{align*}
       <&D_{2N+1}(g)\psi,\phi>_{L^2}\\
        =&<C^{spin}(g,-\frac{2N+n}{2})\circ P_{2N}^{\mathcal{S}(M)}(g)
          \circ D^{spin}(g,\frac{2N+1-n}{2})\psi,\phi>_{L^2}\\
       =&<\psi,i^pi^pC^{spin}(g,-\frac{2N+n}{2})\circ P_{2N}^{\mathcal{S}(M)}(g)
          \circ D^{spin}(g,\frac{2N+1-n}{2})\phi>_{L^2}\\
       =&<\psi,(-1)^p D_{2N+1}(g)\phi>_{L^2},
     \end{align*}
     which completes the proof.
  \end{bew}

  Now we are going to introduce a new family of conformally covariant differential operators acting 
  on sections of the spin tractor bundle. Consider a series of conformally covariant differential operators 
  $D_k(g):\Gamma(S(M,g))\to\Gamma(S(M,g))$, of bi-degree $(\frac{k-n}{2},-\frac{k+n}{2})$, for odd $k\in\N$, 
  not necessarily conformal powers of the Dirac operator. 
  Using these we may define an operator
  \begin{align}
    L_k(g):=\frac{4}{k+1} D^{spin}(g,-\frac{k+n}{2})\circ D_k(g)\circ C^{spin}(g,\frac{k+1-n}{2}),\label{eq:DieL}
  \end{align}
  acting on $\Gamma(\mathcal{S}(M)_g)$. It satisfies the following:
  \begin{theorem}\label{ConfCovOfL}
    For any odd $k\in\N$, the operator $L_k(g)$ is conformally 
    covariant of bi-degree $(\frac{k+1-n}{2},-\frac{k+1+n}{2})$, i.e., 
    for any $\hat{g}=e^{2\sigma}g$ we have that 
    \begin{align*}
      L_k(\hat{g})\left(e^{\frac{k+1-n}{2}\sigma} T^{\mathcal{S}(M)}(g,\sigma)\right)
       =e^{-\frac{k+1+n}{2}\sigma}T^{\mathcal{S}(M)}(g,\sigma)\circ L_k(g).
    \end{align*} 
  \end{theorem}
  The case $k=1$ and $D_1(g)=\slashed{D}$ was found in a joint work with Andreas Juhl 
  analyzing the conformal transformation law for the operator $P_{2}^{\mathcal{S}(M)}(g)$ in detail. 
  \begin{bew}
    This is a direct consequence of its definition \eqref{eq:DieL}. 
  \end{bew}
  
  \begin{bem}
    Note that both operators $P_{2N}^{\mathcal{S}(M)}(g)$ and $L_{2N-1}(g)$ have the same conformal weights,
    see Theorems \ref{ConfCovOfP} and \ref{ConfCovOfL}. 
    Their construction, given in equations \eqref{eq:DieP} and \eqref{eq:DieL}, 
    can be illustrated, for the case $N=2$, as follows:
    \begin{align*}
      \begin{matrix}
        \Gamma(\mathcal{S}^1(M))&\overset{D^0(g,\frac{4-n}{2})}{\xleftarrow{\hspace*{1.5cm}}}
           &\Gamma(\mathcal{S}(M)_g)&
           &\overset{C^{spin}(g,\frac{4-n}{2})}{\xrightarrow{\hspace*{1.5cm}}}&\Gamma(S(M,g))\,\\
        \Box^\nabla_{\frac{2-n}{2}}\downarrow& & 
         \text{\rotatebox[origin=c]{-90}{$\dashrightarrow$}} & & & \downarrow D_3(g)  \\
         \Gamma(\mathcal{S}^1(M))&\underset{C^0(g,-\frac{2+n}{2})}{\xrightarrow{\hspace*{1.5cm}}}
           &\Gamma(\mathcal{S}(M)_g)&
           &\underset{D^{spin}(g,-\frac{n+3}{2})}{\xleftarrow{\hspace*{1.5cm}}}&\Gamma(S(M,g)).
     \end{matrix}
    \end{align*}
    Thus the operators $P^{\mathcal{S}(M)}_4(g)$ and $L_3(g)$ (up to a constant) 
    arise by the dashed arrow depending on the path taken through the diagram. Note, that in general 
    a translation of $L_{2N-1}(g)$ to the spinor bundle vanishes identically, due to 
    $C^{spin}(g,\frac{2N-n}{2})\circ D^{spin}(g,\frac{2N+1-n}{2})=0$, whereas a 
    translation of $P^{\mathcal{S}(M)}_{2N}(g)$ to the spinor bundle 
    yields a conformal power of the Dirac operator. 
  \end{bem}

  Now consider the conformal powers of the Dirac operator $\mathcal{D}_k$, for $k=1,3$, 
  and denote by $\mathcal{L}_k(g)$ the induced conformally covariant operator acting on the spin tractor bundle, 
  given by equation \eqref{eq:DieL}. 
  \begin{theorem}\label{DiffOfLandP}
   On the spin tractor bundle one has that 
   \begin{align*}
     \mathcal{L}_1(g)-P_2^{\mathcal{S}(M)}(g)=&\begin{pmatrix} 0&0\\ \mathcal{D}_3&0\end{pmatrix},\\
     4\mathcal{L}_3(g)-\frac{4}{4-n}P_4(g)=&\begin{pmatrix} 0&0\\ \mathcal{D}_5&0\end{pmatrix},
   \end{align*}
   where in the second difference we have chosen the main part of $P^{\mathcal{S}(M)}_4(g)=P_4(g)+R(g)$.
  \end{theorem}
  The proof based on explicit formulas of the involved operators and can be found in 
  \cite[Theorem $5.48$, Theorem $5.49$]{Fischmann}.

  \begin{bem}
    Theorem \ref{DiffOfLandP} gives also a construction of a conformal third and fifth power of the Dirac 
    operator. It differs to the construction \eqref{eq:ConfDirac}, since we are looking at certain differences of 
    $\mathcal{L}_{2N-1}(g)$ and $P^{\mathcal{S}(M)}_{2N}(g)$, for $N=1,2$, instead of translating 
    $P^{\mathcal{S}(M)}_{2N}(g)$, for $N=1,2$, to the spinor bundle. Of course, translating those 
    differences to the spinor bundle will give us nothing new, since $\mathcal{L}_{2N-1}(g)$ is 
    canceled by translations.  
  \end{bem}

  Now, we come to the polynomial structure of the first examples of conformal powers of the Dirac operator. 
  Using the explicit formulas for $\mathcal{D}_k$, for $k=1,3,5$, we can define 
  differential operators $M_k$, for $k=1,3,5$, by
  \begin{align*}
    M_1:=&\mathcal{D}_1-0\\
    =&\frac 12(g,\nabla^{S(M,g)})+\frac 12(\nabla^{S(M,g)},g\cdot),\\
    M_3:=&\mathcal{D}_3-\mathcal{D}_1^3\\
    =&-(P,\nabla^{S(M,g)})-(\nabla^{S(M,g)},P\cdot),\\
    M_5:=&\mathcal{D}_5-\mathcal{D}_1\mathcal{D}_3\mathcal{D}_1
     -2(\mathcal{D}_1^2\mathcal{D}_3+\mathcal{D}_3\mathcal{D}_1^2)+4\mathcal{D}_1^5\\
    =&4(2P^2+\frac{1}{n-4}B,\nabla^{S(M,g)})+4(\nabla^{S(M,g)},2P^2\cdot+\frac{1}{n-4}B\cdot)\\
      &-2(C,P\cdot)-2(P,C\cdot).
  \end{align*}
  By definition they are first order operators. Just as for each $\mathcal{D}_k$, the $M_k$, for 
  $k=1,3,5$, are formally (anti-) self-adjoint with respect to the $L^2-$scalar product. More interesting, 
  however, is the following result: 
  \begin{theorem}
    On a spin manifold $(M,g)$ of dimension $\neq 4$ we have 
    \begin{align*}
      \mathcal{D}_1=&M_1,\\
      \mathcal{D}_3=&M_1^3+M_3,\\
      \mathcal{D}_5=&M_1^5+M_1M_3M_1+2(M_1^2M_3+M_3M_1^2)+M_5.
    \end{align*}
  \end{theorem}
  
  This structure for the conformal powers of the Dirac operator is very similar to that for the conformal powers 
  of the Laplacian discovered by A. Juhl. He presented a complete series of second order differential operators, 
  such that the GJMS-operators can be written as a polynomial in these operators, see \cite[Theorem $1.1$]{Juhl1}. 
  That series was rediscovered by Fefferman and Graham in \cite{FG2}. 
  
  We believe that there is a completely analogous picture for the conformal powers of the Dirac operator. 
  Hence, 
  it is natural to ask about the nature of $M_k$, for $k\in 2\N-1$. For example, is 
  there a generating function for the series of $M_k$, and how can one understand 
  the coefficients arising in the polynomial description of $\mathcal{D}_k$?

\bibliographystyle{amsalpha}
\bibliography{bibliography}

\end{document}